\documentstyle[12pt]{article}
\newcommand{\sect}[1]{\section{#1}\setcounter{equation}{0}}

\font\mbn=msbm10 scaled \magstep1
\font\mbs=msbm7 scaled \magstep1
\font\mbss=msbm5 scaled \magstep1
\newfam\mbff
\textfont\mbff=\mbn
\scriptfont\mbff=\mbs
\scriptscriptfont\mbff=\mbss
\def\mbf{\fam\mbff}
\def\Re{{\mbf R}}

\def\Z{{\mbf Z}}
\def\Co{{\mbf C}}

\def\P{{\mbf P}}

\def\Di{{\mbf D}}
\newtheorem{Th}{Theorem}[section]
\newtheorem{Lm}[Th]{Lemma}
\newtheorem{C}[Th]{Corollary}

\newtheorem{Proposition}[Th]{Proposition}
\newtheorem{R}[Th]{Remark}
\newtheorem{Ex}[Th]{Example}
\author{Alexander Brudnyi\thanks
{Research supported in part by NSERC.
  }}
\title{LIMIT CYCLES AND THE 
DISTRIBUTION OF ZEROS OF FAMILIES OF ANALYTIC FUNCTIONS}
\date{February 18, 1999}
\begin{document}
\maketitle
\begin{abstract}
{We estimate the expected number of limit cycles situated in a neighbourhood 
of the origin of a planar polynomial vector field.
Our main tool is a distributional inequality for the number of zeros of some 
families of univariate holomorphic functions depending analytically on a 
parameter. We obtain this inequality by methods of Pluripotential Theory. 
This inequality also implies versions of a strong law of large 
numbers and the central limit theorem for a probabilistic 
scheme associated with the distribution of zeros.  
} 
\end{abstract}
\sect{\hspace*{-1em}. Introduction.}
{\bf 1.1.} The second part of Hilbert's sixteenth problem 
asks whether the number of isolated closed trajectories
(limit cycles) of a planar polynomial vector field is always
bounded in terms of its degree. This is the most prominent 
finiteness problem, related to a fairly general class of algebraic
differential equations and is one of the few Hilbert's problems, 
which remain unsolved. Recent results of Y. Il'yashenko [I] and
of \'{E}calle, Martinet, Moussu and Ramis [EMMR] give global finiteness
of the number of limit cycles for each individual vector field (but leave
open the question of the existence of a bound depending on the degree only).
In our paper we consider the local version of the problem in which
one asks for explicit bounds (in terms of degree only) on the number
of limit cycles  situated in a small neighbourhood of a singular point of
the vector field. Except for the result of Bautin [B], the answer to this 
problem is not known. According to Smale (see [S]) the global estimate
should be polynomial in the degree $d$ of the components of the vector field.
Our estimate below (local case) is ``in the mean'' but gives a substantially
better estimate of the order $(\log d)^{2}$. Moreover, Example \ref{lef}
shows that locally the maximal number of limit cycles can be essentially bigger
than $(\log d)^{2}$. We now formulate our result precisely.

Consider a system of ODE's in $\Re^{2}$
\begin{equation}\label{poi}
\begin{array}{c}
\dot x=-y+F(x,y)\\
\dot y=\ x+G(x,y)
\end{array}
\end{equation}
where $F$ and $G$ are polynomials of degree $d$ whose Taylor expansions at 0 
begin with terms of degree $\geq 1$. The origin $(0,0)$ is a singular point 
of (\ref{poi}). 
Assume that 
$$
F=\sum_{k=1}^{d}F_{k},\ \ \ \ \ G=\sum_{k=1}^{d}G_{k}
$$
where
$$
F_{k}(x,y)=\sum_{i=0}^{k}a_{ki}x^{i}y^{k-i},\ \ \ \ \ 
G_{k}(x,y)=\sum_{i=0}^{k}b_{ki}x^{i}y^{k-i}
$$
with real $a_{ki},b_{ki}$. 
Let
$$
|F_{k}|:=\left(\sum_{i=0}^{k}|a_{ki}|^{2}\right)^{1/2},\ \ \ \ \
|G_{k}|:=\left(\sum_{i=0}^{k}|b_{ki}|^{2}\right)^{1/2}\ .
$$
Assume also that
\begin{equation}\label{n}
\sum_{k=1}^{d}(a^{k-1}|F_{k}|)^{2}+
\sum_{k=1}^{d}(a^{k-1}|G_{k}|)^{2}\leq N^{2}\ .
\end{equation}
Condition (\ref{n}) determines an ellipsoid $E(a,N)\subset\Re^{s}$,
$s:=d(d+3)$, in the space of parameters. In what follows we 
identify pairs $F,G$ with points 
$v\in\Re^{s}$. Further,
for any $F,G$ corresponding to $v\in\Re^{s}$ denote by $C(v,K)$ 
the number of limit cycles
of (\ref{poi}) in the disk $D_{K}:=\{(x,y)\in\Re^{2};\
|x|^{2}+|y|^{2}\leq K^{2}\}$. 
We will write $C(v,K)=+\infty$ if any trajectory in $D_{K}$ is closed
(the case of the center). Let $|\cdot |$ denote the
Lebesgue measure on $\Re^{s}$.\\
{\bf Theorem A}\ \
{\em Let $N\leq\frac{1}{192\pi d^{2}}$ be a positive number. 
There are absolute positive constants $c_{1},c_{2}$ such that for any
$T\geq 0$ } 
$$
|\{v\in E(a,N);\ C(v,a/2)\geq T\}|\leq c_{1}se^{-c_{2}T/\log d}
\cdot |E(a,N)| \ .
$$
\begin{C}\label{average}
Under assumptions of Theorem A the expected number of limit cycles
of a random vector field (\ref{poi}) in the disk $D_{a/2}$ is bounded from 
above by $c(\log d)^{2}$ with an absolute constant $c>0$.
\end{C}
The corollary follows from the distributional inequality of Theorem A 
which we apply to estimate the integral \ \ 
$[\int_{E(a,N)}C(v,a/2)dv]/|E(a,N)|$ by $c(\log d)^{2}$.
\begin{Ex}\label{lef}(The maximal number of limit cycles in 
$D_{a/2}$ can be essentially bigger than $(\log d)^{2}$.)
{\rm Consider the system}
$$
\begin{array}{c}
\dot x=-y+xf(x^{2}+y^{2})\\
\dot y=\ x+yf(x^{2}+y^{2})
\end{array}
$$
{\rm where $f$ is a real polynomial whose degree $l$ is 
the integer part of $(d-1)/2$. Assume also that $f$ has
$l$ different positive roots $x_{1},...,x_{l}$ in open interval 
$(0,1/2)$. Then the above system has $l$ limit cycles
$s_{i}:=\{(x,y)\in\Re^{2};\ |x|^{2}+|y|^{2}=x_{i}\}$ for
$i=1,...,l$ in $D_{1/2}$ (see, e.g. [Le], Ch.X, Sec.5).
Clearly we can choose coefficients of $f$ so small that 
the vector of coefficients of polynomials determining the above 
vector field belongs to $E(1,N)$.}
\end{Ex}

{\bf 1.2.} The proof of Theorem A is based on a distributional inequality 
for the number of zeros of some families of univariate holomorphic functions
depending analytically on parameter. Let us note 
that in recent years the problem of estimating the number of zeros 
for families of analytic functions depending analytically on a parameter
has been extensively studied in connection with different aspects of modern 
analysis; e.g. 2nd part of Hilbert's 16th problem, dynamical and control 
systems etc (see  [BY], [FY], [IY],[G],[NY], [Y]). 
However, in many cases the
existing estimates are essentially differ from desired. 
In view of the inequality of Theorem B below it seems likely
that ``good'' estimates in similar problems can be obtained
in the mean. To formulate the result,
let $f:=\{f_{v}\ ; v\in B_{c}(0,r)\}$, $r>1$, be a family of holomorphic
in the open unit disk $\Di_{1}$ functions depending holomorphically on 
parameter $v$ varying in the open Euclidean ball 
$B_{c}(0,r)\subset\Co^{N}$. Assume that for some \ 
$\Di_{s}:=\{z\in\Co\ ; |z|<s\}$ with $s<1$
\begin{equation}\label{1} 
\displaystyle
\sup_{v\in B_{c}(0,r)}\sup_{z\in\Di_{1}}|f_{v}(z)|\leq M<\infty\ \ \ 
{\rm and}\ 
\ \ \sup_{v\in B_{c}(0,1)}\sup_{z\in\Di_{s}}|f_{v}(z)|\geq 1 . 
\end{equation}
Denote the set of these functions by ${\cal H}(M,r,s)$.  

Set, further,  for $f\in {\cal H}(M,r,s)$
\begin{equation}\label{2}
{\cal N}_{f}(v):=\#\{z\in\overline{\Di}_{s};\ \ f_{v}(z)=0\};
\end{equation}
in addition, ${\cal N}_{f}(v)=+\infty$, if $f_{v}=0$ identically on 
$\Di_{1}$. Let $|\cdot |$ denote the Lebesgue measure on $\Co^{N}$.\\
{\bf Theorem B}\ \
{\em For every $T\geq 0$ inequality
$$
|\{v\in B_{c}(0,1) \ ;\ {\cal N}_{f}(v)\geq T\}|\leq 
N c_{1}e^{-c_{2}T/log M}\cdot |B_{c}(0,1)|,
$$
holds with constants $c_{1}, c_{2}$ depending only on $s,r$.}
\begin{R}\label{fef}
According to the doubling inequality from [FN] and 
inequality (\ref{number}) below, the function ${\cal N}_{f}$ is uniformly 
bounded on the set $B_{c}(0,1)\setminus V_{f}$, where 
$V_{f}:=\{v\in B_{c}(0,1)\ ; f_{v}\equiv 0\}$. As we will show 
(see Proposition \ref{pr1} below) the set  $V_{f}$ has measure 0. 
The following example shows, nevertheless, 
that ${\cal N}_{f}$ can assign arbitrary big values on 
$B_{c}(0,1)\setminus V_{f}$.
\end{R} 
\begin{Ex}\label{ex1} {\rm Let 
${\cal O}^{N}(\Di_{1};B_{c}(0,1))$ denote the set of holomorphic 
mappings $f:\Di_{1}\longrightarrow B_{c}(0,1)\subset\Co^{N}$.  
For any  $f\in {\cal 
O}^{N}(\Di_{1};B_{c}(0,1))$ consider the number of intersections of 
$f(\overline{\Di}_{s})$, $s<1$,
with hyperplane $\{(z_{1},...,z_{N})\in\Co^{N};\ 
l(a,z):=a_{1}z_{1}+...+a_{n}z_{N}+a_{N+1}=0\}$, where
$a=(a_{1},...,a_{N+1})\in\Co^{N+1}$. It
coincides with the number of zeros in $\overline{\Di}_{s}$ of the 
function 
$l(a,f(z))$.  For any $f\in {\cal O}^{N}(\Di_{1};B_{c}(0,1))$ the 
function $l(a,f(z))$ satisfies (\ref{1}) with $M=3, r=2$  
if $a\in B_{c}(0,2)$.  Let $B_{k}(z)$ be a 
Bl\"{a}shke product with $k$ zeros in $\Di_{s}$.   
Then the mapping $f(z):=(B_{k}(z),0,...,0)$ belongs to ${\cal 
O}^{N}(\Di_{1};B_{c}(0,1))$ and $l(1,0,...,0,f(z))$ has $k$ zeros 
in $\Di_{s}$.}
\end{Ex}

We prove Theorem B in Section 2 and then in Section 3 we prove 
Theorem A.\\

{\bf 1.3.} In 1943 M.Kac [Ka] proved that the expected number of real 
zeros of a random polynomial of degree $n$ whose coefficients are 
independent standard normals asymptotically equals 
$\frac{2}{\pi}\log n$ as $n\to\infty$.  This theorem was a starting 
point for many results on zeros of random polynomials and other 
functions {\em linearly} depending on random variables (e.g., see [EK]). 
Nonlinear problems of this kind appear
in many important fields of pure and applied mathematics. Here we
mention only a class of problems related to distribution of 
eigenvalues of random matrices (among other applications see [Gi], [M], [Mu] 
and [EK] for the corresponding results and applications to 
quantum physics and multivariate statistics).

Based on Theorem B we study similar problems for the family of 
analytic functions $\{f_{v}\}$ of Section 1.2.

Let 
$f^{(k)}=\{f_{v}^{(k)}\ ; v\in B_{c}(0,r)\}$, $k=1,2,...$, 
be a sequence of functions from ${\cal H}(M,r,s)$. 
Consider a sequence $\{N_{k}\}_{k=1}^{\infty}$ of random 
variables defined on the probability space $B_{c}(0,1)$ as follows.
For a nonnegative integer $l$ probability 
\begin{equation}\label{distrib}
P(N_{k}=l):=
\frac{1}{|B_{c}(0,1)|}|\{v\in B_{c}(0,1)\ ; {\cal N}_{f^{(k)}}(v)=l\}|.
\end{equation}
Denote, as usual, the expectation of $N_{k}$ by $E(N_{k})$ and the variance
of $N_{k}$ by $D(N_{k})$.
\begin{Th}\label{te2}
There are positive constants $c=c(r,s),\ \widetilde c=\widetilde c(r,s)$ such
that
$$
\sup_{k}E(N_{k})\leq c\log M\log(N+1),\ \ \ \
\sup_{k}D(N_{k})\leq\widetilde c(\log M\log(N+1))^{2}.
$$
\end{Th}

Let us consider $\Omega=\prod_{k\geq 1}^{\infty}\Z_{+}$ and the product
probability $\P$ on $\Omega$ associated with distribution (\ref{distrib}).
We study the following probabilistic scheme related to distribution $\P$.
For every $k=1,2,...$ we decompose $B_{c}(0,1)$ into (finite) number of domains
where ${\cal N}_{f^{(k)}}$ is constant (removing a set of measure zero where
${\cal N}_{f^{(k)}}=+\infty$), see Remark \ref{fef}. Then for each 
$k=1,2,...$ we choose at random and independently one of these domains.  
Theorem \ref{te2} guarantees for this scheme fulfillment of 
\begin{C}\label{col1}
The following inequality
$$
\limsup_{n\to\infty}\frac{1}{n}\sum_{k=1}^{n}N_{k}\leq
c\log M\log(N+1)
$$
holds with $\P$ probability one (see also Example \ref{ex1'}).
\end{C}

Assume in addition that the family 
$\{f^{(k)}\}\subset {\cal H}(M,r,s)$ satisfies the following conditions:\\
there are a constant $\delta>0$ and an open disk $\Di_{s'}$ with $s'<s$ such 
that for every $k$ 
$$ 
\begin{array}{lr}
\displaystyle
(a)\ \  
\min_{z\in\overline{\Di}_{s}}|f_{v}^{(k)}(z)|>\delta\\ 
{\rm for\ some}\ \ v=v(k)\in\overline{B}_{c}(0,1);\\ 
\displaystyle
(b)\ \ 
\max_{z\in\Di_{s'}}|f_{v'}^{(k)}(z)|>\delta \ \ {\rm and}\ \ 
f_{v'}^{(k)}(z)=0\\
{\rm for\ some}\ \  v'=v'(k)\in\overline{B}_{c}(0,1)\ \ {\rm and}\ \ 
z=z(k)\in\Di_{s'}. 
\end{array}
$$
Under these assumptions and in the above notations the following result holds.
\begin{Th}\label{central} 
The sequence $\{N_{k}\}_{k=1}^{\infty}$ 
of independent random variables satisfies the normal distribution law, 
that is, 
$$ 
P\{\frac{1}{B_{n}}\sum_{k=1}^{n}(N_{k}-E(N_{k}))<x\}\longrightarrow
\frac{1}{\sqrt{2\pi}}\int_{-\infty}^{x}e^{-z^{2}/2}dz,\ \ \ (n\to\infty)
$$
uniformly in $x$. 
Here $B_{n}^{2}:=\sum_{k=1}^{n}D(N_{k})$ is the sum of the variances of 
$N_{k}$.
\end{Th}
\begin{Ex}\label{ex2}
{\rm  
Let $\{f_{k}(z):=(cz,f_{1,k}(z),...,f_{N-1,k}(z))\}_{k=1}^{\infty}$,
$c\neq 0$, be
a sequence of holomorphic mappings from ${\cal O}^{N}(\Di_{1};B_{c}(0,1))$,
i.e. $\sup_{z\in\Di_{1}}(|cz|^{2}+\sum_{i=1}^{N-1}|f_{i,k}(z)|^{2})<1$ for
each $k\geq 1$ (see Example \ref{ex1}). We set 
$f_{v}^{(k)}(z):=l(v, f_{k}(z))$.
Then, clearly, the sequence $\{f^{(k)}\}_{k\geq 1}$ satisfies the
conditions of Theorem \ref{central}.}
\end{Ex}

Consider another theorem of this type generalizing a classical result 
on behavior of zeros of random polynomials.

Let $P_{k,v}(z)=\sum_{i=0}^{k}a_{ik}(v)z^{i}$ be a polynomial in 
$z\in\Co$ and $v\in\Co^{n(k)}$ of degree $d(k)$ in $v$ satisfying
$$
\begin{array}{l}
\displaystyle
(a)\ \ d(k)\log n(k)=o(k)\ \ \ {\rm as}\ \ \ k\to\infty;\\
\displaystyle
(b)\ \ \sup_{||v||\leq 1}|a_{ik}(v)|\leq 1,\ \ \ {\rm if}\ \ \ 1\leq 
i\leq k-1;\\
\displaystyle
(c)\ \ \sup_{||v||\leq 1}|a_{ik}(v)|=1,\ \ \ {\rm if}\ \ \ i=0,k.
\end{array}
$$
Let $\epsilon >0$ be arbitrary small. Denote by  $A_{\epsilon}$ the 
annulus $\{z\in\Co\ ;\ 1-\epsilon <|z|<1+\epsilon\}$. Consider the 
counting function $\widetilde {\cal N}_{k}(v)$ defined by
$$
\widetilde {\cal N}_{k}(v):=\#\{z\in A_{\epsilon}\ ; P_{k,v}(z)=0\}.
$$
Further, determine as before  a random variable $\widetilde N_{k}$ 
defined on the probability space $B_{c}^{n(k)}:=\{v\in\Co^{n(k)}\ ; 
||v||<1\}$ by
$$
P(\widetilde N_{k}=l):=\frac{1}{|B_{c}^{n(k)}|}mes_{2n(k)}(\{v\in 
B_{c}^{n(k)}\ ; \widetilde {\cal N}_{k}(v)=l\}),
$$
where $l=0,1,...$.
\begin{Th}\label{polin}
Expectation $E(\widetilde N_{k})$ of the numbers of zeros of $P_{k,v}$ 
in $A_{\epsilon}$ is bounded from below by $k(1-o_{\epsilon}(1))$ as 
$k\to\infty$.
\end{Th}
\begin{R}\label{circle}
Thus, as in the linear case $({\rm \ i.e.,\ for\ } 
a_{ik}(v)=v_{i+1},\ 0\leq i\leq 
n(k)-1:=k)$ the zeros of polynomials $P_{k,v}$ concentrate on the unit 
circle, as their degrees grow. It is possible to prove that if 
$\{a_{ik}\}_{i=0}^{k}$ are homogeneous and mutually independent on 
$B_{c}^{n(k)}$\ $(k=1,2,...)$ then roots of random polynomials 
$P_{k,v}$ are asymptotically uniformly distributed on the unit circle 
as $k\to\infty$.
\end{R} 
\sect{\hspace*{-1em}. Proof of Theorem B.}
Let us verify, first, that the set
$$
\omega_{T}:=\{v\in \overline{B_{c}}(0,1)\ ; {\cal N}_{f}(v)\geq
T\}
$$
is measurable. 
\begin{Proposition}\label{pr1}
For every $f\in {\cal H}(M,r,s)$ there is a closed subset 
$V=V_{f}\subset B_{c}(0,r)$ of measure $0$ such that ${\cal 
N}_{f}$ is upper semicontinuous on $B_{c}(0,r)\setminus V_{f}$.  
\end{Proposition} 
{\bf Proof.}  
Let 
$$ 
V_{f}:=\{v\in B_{c}(0,r);\ \ f_{v}\equiv 0\ on\ \Di_{1}\}.  
$$ 
Prove that $V_{f}$ 
is a proper complex analytic subset of $B_{c}(0,r)$. To this end we set 
$$
\widetilde V_{f}:=\bigcap_{k\geq 0}
\{(v,z)\subset B_{c}(0,r)\times\Di_{1};\ \ 
\frac{\partial ^{k}}{\partial z^{k}} f_{v}(z)=0\}.  
$$
According to the definition ${\widetilde V}_{f}$ is a complex 
analytic subset of $B_{c}(0,1)\times\Di_{1}$.
If $\Pi:B_{c}(0,r)\times\Di_{1}\longrightarrow B_{c}(0,r)$ is the 
natural projection, then $\Pi (\widetilde{V}_{f})$, clearly, coincides 
with $V_{f}$. Hence $\Pi^{-1}(v)\subset\widetilde{V}_{f}$ for every 
$v\in V_{f}$. Therefore $V_{f}$ is biholomorphically isomorphic to the 
set $\widetilde{V}_{f}\cap (B_{c}(0,r)\times\{0\})$ which is 
analytic as intersection of analytic sets. It remains to check that 
$V_{f}$ is a proper subset of $B_{c}(0,r)$. If, to the contrary,  
$\widetilde V_{f}$ coincides with $B_{c}(0,r)\times\Di_{1}$, then 
$\{f_{v}\}_{v}=\{0\}$ which 
contradicts to the second inequality of (\ref{1}). 

Prove now that  ${\cal N}_{f}$ is upper semicontinuous on 
$B_{c}(0,r)\setminus V_{f}$. This means that 
\begin{equation}\label{upper}
\limsup_{k\to\infty}{\cal N}_{f}(v_{k})\leq 
{\cal N}_{f}(v)
\end{equation}
for every $\{v_{k}\}\subset B_{c}(0,r)\setminus V_{f}$ converging to 
$v\in B_{c}(0,r)\setminus V_{f}$.  
Denote by $\gamma_{v}$ the boundary of the disk $\{z;\ |z|\leq r\},\ 
r>s$, containing as the same number of zeros of $f_{v}$ as 
$\overline{\Di}_{s}$. Then we can represent ${\cal N}_{f}$ by 
\begin{equation}\label{integ1} 
\frac{1}{2i\pi}\int_{\gamma_{v}}\frac{\frac{\partial}{\partial z}
f_{v}(z)}{f_{v}(z)}dz.  
\end{equation}
By continuity of $f_{v}(z)$ there is an open connected neighbourhood 
$U_{v}$ of $v$ such that the function $f_{w}$ has no zeroes on 
$\gamma_{v}$ for each $w\in U_{v}$. Therefore the right hand side of
(\ref{integ1}) is well defined for such $f_{w}$ and is a continuous  
function in $w$ on $U_{v}$. Moreover, this function assigns only integer 
values and $U_{v}$  is connected. Hence,
\begin{equation}\label{2.3}
\frac{1}{2i\pi}\int_{\gamma_{v}}\frac{\frac{\partial}{\partial z}
f_{w}(z)}{f_{w}(z)}dz={\cal N}_{f}(v)  
\end{equation}
for all $w\in U_{v}$. Let now $\{v_{k}\}$ be a sequence from 
(\ref{upper}). Because of the precompactness of the family $\{f_{v}\}$
in the topology of uniform convergence on compact subsets of $\Di_{1}$, 
see (\ref{1}), we can assume that $\{f_{v_{k}}\}$ converges to 
$f_{v}$ in this topology. In particular, $\{f_{v_{k}}|_{\gamma_{v}}\}$ 
converges uniformly to $f_{v}|_{\gamma_{v}}$. Let $k_{0}$ be such that 
$v_{k}$ belongs to $U_{v}$ for any $k\geq k_{0}$. Then the number of 
zeros ${\cal N}_{k}$ of $f_{v_{k}},\ k\geq k_{0}$, insides 
$\gamma_{v}$ coincides with ${\cal N}_{f}(v)$. Moreover, 
$\overline{\Di}_{s}$ is also containing in the disk bounded by 
$\gamma_{v}$; therefore
$$
{\cal N}_{f}(v_{k})\leq {\cal N}_{k}={\cal N}_{f}(v).
$$

The proposition is proved.\ \ \ \ \ $\Box$
\begin{R}\label{open}
Let 
$$
{\cal \widetilde N}_{k}(v):=\{z\in\Di_{s}\ ; f_{v}(z)=0\}
$$
denote the number of zeros in the {\em open} disk\ $\Di_{s}$. Then
by similar arguments one can establish that 
${\cal \widetilde N}_{k}$ is lower semicontinuous on the set
$B_{c}(0,r)\setminus V_{f}$.
\end{R}

Verify now that $\omega_{T}$ is a compact subset of 
$\overline{B}_{c}(0,1)$. In fact, let 
$\{v_{k}\}\subset\omega_{T}$ be a sequence converging to 
$v\in\overline{B}_{c}(0,1)$. If $v\in B_{c}(0,r)\setminus V_{f}$ then,
according to upper semicontinuity of ${\cal N}_{f}$, it also belongs
to $\omega_{T}$. Otherwise, $v\in V_{f}$ and therefore
${\cal N}_{f}(v)=+\infty$ implying $v\in\omega_{T}$. 

The next part of the proof involves arguments of Pluripotential Theory. 
We, first, recall the following classical estimate for the number of 
zeros of a holomorphic on $\Di_{1}$ function $h$ in the closed disk 
$\overline{\Di}_{s}$ (a consequence of Jensen's inequality).
 
Let $M_{1}:=\sup_{\Di_{\frac{s+1}{2}}}\log|h|$ and 
$M_{2}:=\sup_{\Di_{s}}\log|h|$.  Then 
\begin{equation}\label{number} 
\#\{z\in \overline{\Di}_{s};\ \ h(z)=0\}\leq c(M_{1}-M_{2}) 
\end{equation} 
with $c=c(s)$.

Let now the family $f=\{f_{v}\}$ belongs to ${\cal H}(M,r,s)$. 
For a fixed $v\in B_{c}(0,r)$ we set 
$$ 
\displaystyle 
M_{f,1}(v):=\frac{\sup_{\Di_{\frac{s+1}{2}}}\log|f_{v}|-\log M}{\log M}
\ \ \ and\ \ \ 
M_{f,2}(v):=\frac{\sup_{\Di_{s}}\log|f_{v}|-\log M}{\log M}.
$$
Since $f_{v}(z)$ is continuous on $\Di_{1}\times B_{c}(0,r)$, 
the functions $e^{M_{f,1}}$ and $e^{M_{f,2}}$ are continuous on 
$B_{c}(0,r)$. In particular, $M_{f,1}$ and $M_{f,2}$ are upper semicontinuous.
Moreover, each of these functions is, by definition, supremum of a family of 
plurisubharmonic on $B_{c}(0,r)$ functions. Therefore they are also 
plurisubharmonic.
Further, inequality (\ref{number}) gets
\begin{equation}\label{estim}
{\cal N}_{f}(v)\leq c\log M(M_{f,1}(v)-M_{f,2}(v)).
\end{equation}
Note that from (\ref{1}) it follows that
$\sup_{B_{c}(0,r)}M_{f,1}\leq 0$. Moreover,
inequalities (\ref{1}) imply that there is
a point $x_{f}$ in the open ball $B_{c}(0,1)$ such that
$$
\sup_{B_{c}(0,1)}M_{f,2}>M_{f,2}(x_{f})\geq -2.
$$
From Lemma 3 of [BG] it follows that there exists a ray $l_{f}$ with 
the  origin at $x_{f}$ such that 
\begin{equation}\label{BG}
\frac{mes_{1}(B_{c}(0,1)\bigcap l_{f})}
{mes_{1}(\omega_{T}\bigcap l_{f})}\leq
\frac{2N|B_{c}(0,1)|}{|\omega_{T}|}.
\end{equation}
Consider the one-dimensional affine complex line $l_{f}'$ containing
$l_{f}$. 
Let $D_{1}, D_{2}$ be the intersections of $l_{f}'$ with $B_{c}(0,r)$
and $B_{c}(0,1)$, respectively. If $z_{f}$ denotes the point of 
$l_{f}'$ such 
that $d_{f}:=\displaystyle dist(0,l_{f}')=|z_{f}|$, then $D_{1}, D_{2}$ are
open disks in $l_{f}'$ centered at $z_{f}$ with radii 
$$
r_{1}:=\sqrt{r^{2}-d_{f}^{2}},\ \ \ \ 
r_{2}:=\sqrt{1-d_{f}^{2}},
$$
respectively. Since $r>1$, the ratio $r_{1}/r_{2}\geq r$. Consider now 
the disks $$ \widetilde D_{2}:=\frac{1}{rr_{2}}(D_{2}-z_{f})\ \ and\ \ 
\widetilde D_{1}:=\frac{1}{rr_{2}}(D_{1}-z_{f})
$$
We can thought of them as open disks in $\Co$.
Then, clearly, 
$$
\Di_{1/r}=\widetilde D_{2}\subset\Di_{1}\subset\widetilde D_{1}.
$$
Further, consider  the function $M'$ defined on $\Di_{1}$ by
$M'(z):=M_{f,2}(rr_{2}z+z_{f})$. Since $M_{f,2}$ is 
non-positive plurisubharmonic, the function $M'$ has the same 
property and satisfies: 
$$
\sup_{\Di_{1/r}}M'\geq -2.
$$ 
Set 
$K:=\omega_{T}\cap l_{f}$ and $K':=\frac{1}{rr_{2}}(K-z_{f})$ and
$I_{r}=\frac{1}{rr_{2}}((D_{2}\cap l_{f})-z_{f})$.
Then according to Theorem 1.2 of [Br] (see also section 2.3 there) there is a
constant $c(r)>0$ such that
\begin{equation}\label{brudest}
-2\leq\sup_{\Di_{1/r}}M'\leq c(r)\log\frac{4|I_{r}|}{|K'|}+\sup_{K'}M'.
\end{equation}
Noting that 
$$
\sup_{K'}M'=\sup_{K}M_{f,2}\ \ \ {\rm and}\ \ \
\frac{|I_{r}|}{|K'|}=\frac{mes_{1}(B_{c}(0,1)\bigcap l_{f})}
{mes_{1}(\omega_{T}\bigcap l_{f})}
$$
and taking into account (\ref{BG}) we obtain from (\ref{brudest}) 
\begin{equation}\label{middle}
\sup_{K}M_{f,2}\geq -2-c(r)\log\frac{8N |B_{c}(0,1)|}{|\omega_{T}|}.
\end{equation}
Further, assume
without loss of generality that $\sup_{K}M_{f,2}=M_{f,2}(\widetilde 
v)$ for some $\widetilde v\in\omega_{T}\cap l_{f}$. Then
from (\ref{estim}) and (\ref{middle}) it follows that  
$$
\begin{array}{c}
\displaystyle
T\leq {\cal N}_{f}(\widetilde v)\leq c\log 
M(M_{f,1}(\widetilde v)-M_{f,2}(\widetilde v))
\leq -(c\log M)M_{f,2}(\widetilde v)\\
\displaystyle
\leq  c\log M\left(2+c(r)\log\frac{8N |B_{c}(0,1)|}{|\omega_{T}|}\right).
\end{array}
$$
The previous inequality implies
$$
|\omega_{T}|\leq  N c_{1}e^{-c_{2}T/\log M}|B_{c}(0,1)|
$$
with $c_{1},c_{2}$ depending only on $s,r$.

The proof of the theorem is complete.\ \ \ \ \ $\Box$

Let $f\in {\cal H}(M,r,s)$ be such that
\begin{equation}\label{rever}
\sup_{v\in B_{c}(0,1)}\sup_{z\in\Di_{s}}|f_{v}(z)|\geq
\sup_{v\in B(0,1)}\sup_{z\in\Di_{s}}|f_{v}(z)|\geq 1\ 
\end{equation}
where $B(0,1)\subset B_{c}(0,1)\subset\Co^{N}$ is the real Euclidean ball.
Then one can prove
\begin{Th}\label{rea}
For every $T\geq 0$
$$
|\{v\in B(0,1)\ ;\ {\cal N}_{f}(v)\geq T\}|\leq
N\tilde c_{1}e^{-\tilde c_{2}T/log\ M}\cdot |B(0,1)|,\
$$
where $\tilde c_{1},\tilde c_{2}$ are positive constants depending on 
$s$ and $r$ only.\\
Here $|\cdot|$ is the Lebesgue measure on $\Re^{N}$.
\end{Th}
{\bf Proof.} The proof repeats word-for-word our proof of 
Theorem B. We have to use only that the above point $x_{f}$ can be taken
from $B(0,1)$ and also instead of (\ref{BG}) the real version of 
Brudnyi-Ganzburg lemma (see [BG])
$$
\frac{mes_{1}(B(0,1)\bigcap l_{f})}{mes_{1}(\omega_{T}\bigcap l_{f})}\leq
\frac{N|B(0,1)|}{|\omega_{T}|}
\ \ \ \ \ \Box
$$
\sect{\hspace*{-1em}. Proof of Theorem A.}
By the change of variables
$x\mapsto x/a$, $y\mapsto y/a$  we reduce (\ref{poi}) to 
the equivalent system with polynomial terms $F=\sum_{k=1}^{d}F_{k}$,
$G=\sum_{k=1}^{d}G_{k}$ satisfying
$$
\sum_{k=1}^{d}|F_{k}|^{2}+\sum_{k=1}^{d}|G_{k}|^{2}\leq N^{2}
\ .
$$
Hence it suffices to prove the theorem for $a=1$. Note that in this 
case the ellipsoid $E(a,N)$ coincides with the Euclidean ball
$B(0,N)\subset\Re^{s}$ with center 0 and radius $N$. So we must 
estimate the number of limit cycles $C(v,1/2)$ in the disk  
$D_{1/2}=\{(x,y)\in\Re^{2};\ x^{2}+y^{2}\leq 1/4\}$.

Writing system (\ref{poi}) in polar coordinates $x=r\cos\theta$, 
$y=r\sin\theta$ we get
$$ 
\begin{array}{c} 
\dot x=\dot r\cos\theta-r\sin\theta\cdot\dot \theta=-r\sin\theta+ 
F(r\cos\theta,r\sin\theta),\\
\dot y=\dot r\sin\theta+r\cos\theta\cdot\dot \theta=\ r\cos\theta+
G(r\cos\theta, r\sin\theta).
\end{array}
$$
Multiplying the first equation by $\cos\theta$, the second by $\sin\theta$
and adding we get
$$
\begin{array}{c}
\dot r=F(r\cos\theta,r\sin\theta)\cos\theta+G(r\cos\theta,r\sin\theta)
\sin\theta=\\
\displaystyle
\sum_{k=1}^{d}r^{k}(F_{k}(\cos\theta,\sin\theta)\cos\theta+G_{k}(\cos\theta,\sin\theta)\sin\theta)
=\sum_{k=1}^{d}r^{k}f_{k}(\theta)\ .
\end{array}
$$
Similarly,
$$
\dot \theta=
1+\sum_{k=1}^{d}r^{k-1}
(-F_{k}(\cos\theta,\sin\theta)\sin\theta+G_{k}(\cos\theta,\sin\theta)
\cos\theta)=1+\sum_{k=1}^{d}r^{k-1}g_{k}(\theta)\ .
$$
Finally, we get
\begin{equation}\label{zam}
\frac{dr}{d\theta}=\frac{\dot r}{\dot \theta}=
\frac{r P(r,\theta)}{1+Q(r,\theta)}
\end{equation}
where
$$
P(r,\theta):=\sum_{k=1}^{d}r^{k-1}f_{k}(\theta)\ \ \ \ \
Q(r,\theta):=\sum_{k=1}^{d}r^{k-1}g_{k}(\theta)\ .
$$
The functions $P$ and $Q$ is clear to depend on $v=(F,G)\in\Re^{s}$
linearly. So we will write $f_{k}(\theta):=f_{k}(v,\theta)$,  
$g_{k}(\theta):=g_{k}(v,\theta)$, $k=1,...,d$, and $P(r,\theta):=
P(v,r,\theta)$, $Q(r,\theta):=Q(v,r,\theta)$, $v\in\Re^{s}$.
\begin{Lm}\label{period}
Let $v\in B(0,N)$, $N\leq\frac{1}{192\pi d^{2}}$. Then the trajectory 
of system (\ref{poi}) in the disk $D_{1}$ is closed if and only if
the corresponding solution of (\ref{zam}) satisfies $r(2\pi)=r(0)$
(and hence is periodic).
\end{Lm}
{\bf Proof.} It suffices to check that denominator $1+Q(v,r,\theta)$
differs from 0 if $v\in B(0,N)$, $0\leq r\leq 1$, $\theta\in\Re$. By 
definition, 
$$ 
G_{k}(\cos t,\sin t)=\sum_{i=0}^{k}b_{ki}(\cos t)^{i}(\sin t)^{k-i}
$$
with $\sum_{i=0}^{k}|b_{ki}|^{2}\leq |v|^{2}$, where $|\cdot |$ denotes
the Euclidean norm in $\Re^{s}$. Hence for $\theta\in\Re$
$$
|G_{k}(\cos\theta,\sin\theta)|
\leq\sqrt{\sum_{i=0}^{k}\frac{|b_{ki}|^{2}}{C_{k}^{i}}} 
\cdot\sqrt{\sum_{i=0}^{k}C_{k}^{i}(\cos\theta)^{2i}
(\sin\theta)^{2(k-i)}}
$$
and therefore
\begin{equation}\label{Gk}
\max_{\theta}|G_{k}(\cos\theta,\sin\theta)|\leq |v|,\ \ \
i=0,1,...,k\ .
\end{equation}
Similarly we get
\begin{equation}\label{Fk}
\max_{\theta}|F_{k}(\cos\theta,\sin\theta)|\leq |v|,\ \ \
i=0,1,...,k\ .
\end{equation}
From (\ref{Gk}) we now have for $v\in B(0,N)$
$$
\max_{0\leq r\leq 1, 0\leq\theta\leq 
2\pi}|Q(v,r,\theta)|\leq 2d|v|\leq 2dN<1,
$$
since $N\leq\frac{1}{192\pi d^{2}}$.\ \ \ \ \ $\Box$  

Assuming that $v$ belongs to the complex Euclidean ball 
$B_{c}(0,2N)\subset\Co^{s}$ we 
extend $f_{k}$ and $g_{k}$ to 
holomorphic in $v$ functions on $B_{c}(0,2N)$. This also gives an extension
of $P$ and $Q$ in (\ref{zam}) to holomorphic in $v$ functions
on $B_{c}(0,2N)$. Consider the complexification of equation 
(\ref{zam}):
\begin{equation}\label{zam1} 
\frac{dr}{dt}=\frac{rP(v,r,t)}{1+Q(v,r,t)} 
\end{equation} 
where $v\in B_{c}(0,2N)$, $t\in\Re$.
We are looking for a complex valued solution $r$ of (\ref{zam1}) defined 
in interval $I=[0,2\pi]\subset\Re$ and satisfying the initial 
condition $r(0)=w\in\Di_{3/4}$. Solve (\ref{zam1}) by the method of 
successive approximations.  Namely, set $r_{0}(v,w,t)=w$ and let 
$$ 
r_{n+1}(v,w,t)\!:=\!w+\int_{0}^{t}\!\frac{r_{n}(v,w,\theta)
P(v,r_{n}(v,w,\theta),\theta)}{1+Q(v,r_{n}(v,w,\theta),\theta)}d\theta
\ \ \ (v\in B_{c}(0,2N),w\in\Di_{3/4}, t\in I).
$$
For a function $f$ defined in $B_{c}(0,2N)\times\Di_{h}\times I$, 
$h>0$, we set 
$$ 
||f||_{h}:=\sup_{v\in B_{c}(0,2N), 
z\in\Di_{h},\theta\in I}|f(v,z,\theta)|\ .  
$$ 
\begin{Proposition}\label{approxim}
The sequence $\{r_{n}\}_{n\geq 0}$ converges uniformly on 
$B_{c}(0,2N)\times\Di_{3/4}\times I$ to a function $r$ which is
analytic on this set. Function $r$ uniquely solves (\ref{zam1}) 
with the initial condition $r(v,w,0)=w$ and satisfies\\
(a)\ \ \ $||r||_{3/4}\leq 1$;\ \ (b)\ \ \ 
$r(\cdot,\cdot, t)$ is holomorphic in $B_{c}(0,2N)\times\Di_{3/4}$ for 
any $t\in I$.  
\end{Proposition} 
{\bf Proof.} We begin with 
\begin{Lm}\label{pq}
$$
||P||_{1}\leq 4Nd\leq 1/(48\pi d),\ \ \ \ \
||Q||_{1}\leq 4Nd<1/2 \ .
$$
\end{Lm}
{\bf Proof.} Note that inequalities (\ref{Gk}), (\ref{Fk}) are also
valid for $v\in\Co^{s}$. Applying these inequalities with
$v\in B_{c}(0,2N)$ and using definitions of $f_{k}$ and $g_{k}$
we then have
$$
\sup_{t\in\Re}|f_{k}(v,t)|\leq 4N,\ \ \ \ \ 
\sup_{t\in\Re}|g_{k}(v,t)|\leq 4N\ .
$$
Hence
$$
\begin{array}{l}
\displaystyle
||P||_{1}\leq 4N\sum_{i=1}^{d}1^{i}=4Nd\leq 1/(48\pi d),\\
\displaystyle
||Q||_{1}\leq 4N\sum_{i=1}^{d}1^{i}=4Nd<1/2 \ .
\end{array}
$$
by the choice of $N$. \ \ \ \ \ $\Box$

Check now that
\begin{equation}\label{tri}
||r_{n}||_{3/4}<1\ \ \ \ {\rm for\ any}\ n\ .
\end{equation}
Indeed, from inequalities of Lemma \ref{pq} it follows that
$$
8\pi ||P||_{1}+||Q||_{1}<1\ .
$$
Therefore for $v\in B_{c}(0,2N)$, $w\in\Di_{3/4}$, $t\in I$
$$
|r_{n+1}(v,w,t)|=\left|w+\int_{0}^{t}\frac{r_{n}(v,w,\theta)
P(v,r_{n}(v,w,\theta),\theta)}{1+Q(v,r_{n}(v,w,\theta),\theta)}d\theta
\right|\leq\frac{3}{4}+\frac{2\pi ||P||_{1}}{1-||Q||_{1}}<1,
$$
provided $||r_{n}||_{3/4}<1$. Since $|r_{0}|=|w|\leq 3/4$ the 
inequality (\ref{tri}) is proved by induction.

Using (\ref{tri}) we now prove that 
$$
||r_{n+1}-r_{n}||_{3/4}<\frac{1}{2}||r_{n}-r_{n-1}||_{3/4}\ .
$$
Applying the mean-valued inequality we obtain
\begin{equation}\label{difer}
||r_{n+1}-r_{n}||_{3/4}\leq
2\pi ||r_{n}-r_{n-1}||_{3/4}\left|\left|
\frac{\frac{\partial (rP)}{\partial r}(1+Q)
-\frac{\partial Q}{\partial r}(rP)}{(1+Q)^{2}}\right|\right|_{1}\ .
\end{equation}
Since $|r|\leq 1$, we can use the 
classical Bernstein inequality for holomorphic polynomials (in $r$)
$$
\begin{array}{c}
\displaystyle
\left|\frac{\partial}{\partial r}(rP(v,r,\theta))\right|
\leq d\ ||P||_{1},\ \ \ 
\left|\frac{\partial }{\partial r}(Q(v,r,\theta))\right|
\leq (d-1)\ ||Q||_{1}\\ 
\displaystyle
((v,r,\theta)\in B_{c}(0,2N)\times\Di_{1}\times I)
\end{array}
$$
to estimate the right-hand side of (\ref{difer}).
From here and Lemma \ref{pq} we have
$$
||r_{n+1}-r_{n}||_{3/4}\leq ||r_{n}-r_{n-1}||_{3/4}
\frac{4\pi d||P||_{1}(1+||Q||_{1})} 
{(1-||Q||_{1})^{2}}<\frac{1}{2}||r_{n}-r_{n-1}||_{3/4}\ .
$$
Hence the sequence
$\{r_{n}\}_{n\geq 0}$ converges uniformly on 
$B_{c}(0,2N)\times\Di_{3/4}\times I$ to a complex valued 
analytic function $r$ such that 
$||r||_{3/4}\leq 1$. This function
uniquely solves equation (\ref{zam1}) with $r(v,w,0)=w$. Since each
$r_{n}(\cdot,\cdot,t)$ is holomorphic on $B_{c}(0,2N)\times\Di_{3/4}$ for any
$t\in I$,
the function $r(\cdot,\cdot,t)$ is holomorphic on
$B_{c}(0,2N)\times\Di_{3/4}$ for any $t\in I$, as well. 

The proposition is proved.\ \ \ \ \ $\Box$
\begin{R}\label{re}
It is worth noting that restriction of $r$ to
$B(0,2N)\times [0,3/4]\times I$ is the real valued
solution of (\ref{zam}). 
\end{R}

We proceed to the proof of Theorem A for $a=1$.

Let $r=r(v,w,t)$, $v\in B_{c}(0,2N), w\in\Di_{3/4}, t\in I$
be the solution of (\ref{zam1}). Let us rewrite this equation as
\begin{equation}\label{integr}
p(v,w):=r(v,w,2\pi)-w=\int_{0}^{2\pi}\frac{r(v,w,\theta)
P(v,r(v,w,\theta),\theta)}{1+Q(v,r(v,w,\theta),\theta)}d\theta
\end{equation}
To estimate the number of limit cycles $C(v,1/2)$ of 
(\ref{poi}) we must, according to Lemma \ref{period} and Remark 
\ref{re}, estimate the number of zeros of $p(v,\cdot)$ in 
$\Di_{1/2}$.  
From (\ref{integr}) and Lemma \ref{pq} one has 
$$ 
||p||_{3/4}=\sup_{B_{c}(0,2N)\times\Di_{3/4}}|p|\leq 16\pi d N\ .  
$$ 
Consider now the system
$$
\begin{array}{c}
\dot x=-y+(N/2)x\\
\dot y=\ x+(N/2)y
\end{array}
$$
and denote by $v_{0}\in B(0,N)\subset\Re^{s}$ the vector of coefficients
of polynomials which determines this system.
Reducing this system to the equation
$$
\frac{dr}{d\theta}=\frac{N}{2}r
$$
we easily see that its solution is
$$
r(v_{0},w,\theta)=e^{N\theta /2}w.
$$
Then 
$$
|p(v_{0},w)|=|1-e^{\pi N}||w|\geq \pi N|w|\ .
$$
In particular, for $w=1/2$ we have
$$
|p(v_{0},1/2)|\geq\frac{\pi N}{2}\ .
$$
Thus 
\begin{equation}\label{ratio1}
\frac{||p||_{3/4}}{|p(v_{0},1/2)|}\leq 32d .
\end{equation}
Further, we set
$$
f_{v}(z):=\frac{p(Nv,3z/4)}{p(v_{0},1/2)}\ \ \ (|z|<1,v\in B_{c}(0,2)).
$$
Now inequality (\ref{ratio1}) shows that 
$$
\sup_{v\in B_{c}(0,2)}\sup_{z\in\Di_{1}}|f_{v}(z)|\leq 32d .
$$
Moreover,
$$
\sup_{v\in B_{c}(0,1)}\sup_{z\in\Di_{2/3}}|f_{v}(z)|\geq
|f_{v_{0}/N}(2/3)|=1 .
$$
Hence $f$ belongs to ${\cal H}(32d,2,2/3)$ and satisfies (\ref{rever}).
So we can apply Theorem \ref{rea} to estimate ${\cal N}_{f}$ in 
$B(0,1)$ which, in turn, coincides with the number of zeros $N_{p}(v)$ 
of $p(v,\cdot)$, $v\in B(0,N)$,  in the disk $\Di_{1/2}$:  
$$ 
\begin{array}{c} 
\displaystyle
|\{v\in B(0,N)\ ;\ C(v,1/2)\geq T\}|\leq 
|\{v\in B(0,N)\ ;\  N_{p}(v)\geq T\}|=\\
\displaystyle
N^{s}|\{v\in B(0,1)\ ;\ {\cal N}_{f}(v)\geq T\}|\leq
d(d+3)\tilde c_{1}e^{-\tilde c_{2}T/\log (32d)} |B(0,N)|,\
\end{array}
$$
where $\tilde c_{1},\tilde c_{2}$ are absolute positive constants.

The proof is complete.\ \ \ \ \ $\Box$\\
{\bf Proof of Corollary \ref{average}.} The proof is similar to the proof of 
Theorem \ref{te2} below and we do not reproduce it.\ \ \ \ \ $\Box$
\sect{\hspace*{-1em}. Proofs.}
In this section we prove results of Section 1.3.
Many of our implications are well known in Probability Theory but the whole
point is that they require in their assumptions the inequality of Theorem B.
\\
{\bf Proof of Theorem \ref{te2}.}
We, first, estimate the expectation of the random variable $N_{k}$, $k\geq 
1$.
By definition
$$
E(N_{k})=\frac{1}{|B_{c}(0,1)|}\sum_{l=0}^{\infty}l\cdot |\{v\in 
B_{c}(0,1)\ ; {\cal N}_{f^{(k)} 
}(v)=l\}|=\frac{1}{|B_{c}(0,1)|}\int_{B_{c}(0,1)}{\cal 
N}_{f^{(k)}}(v)dv.
$$
The latter integral equals
$$
\frac{1}{|B_{c}(0,1)|}\int_{0}^{|B_{c}(0,1)|}{\cal 
N}_{f^{(k)}}^{*}(t)dt,
$$
where $g^{*}$ denotes the nonincreasing rearrangement of a function
$g$. Since $g^{*}$ is, by definition, the right inverse to the 
distribution function $T\mapsto |\{v\in B_{c}(0,1)\ ; |g(v)|\geq T\}|$,
we easily deduce from Theorem B that
$$
{\cal N}_{f^{(k)}}^{*}(t)\leq
\frac{\log M}{c_{2}}\log\left(\frac{c_{1}N|B_{c}(0,1)|}{t}\right).
$$
Putting together this inequality and the previous identity we
prove that
$$
E(N_{k})\leq c\log M\log(N+1),
$$
where $c$ depends only on $r,s$.

Using similar arguments one can show that
\begin{equation}\label{square}
E(N_{k}^{2})\leq c'(\log M\log (N+1))^{2}
\end{equation}
with $c'$ depending only on $r,s$. Then
$$ 
D(N_{k})\leq 2E(N_{k}^{2})+2(E(N_{k}))^{2}\leq\widetilde c
(\log M\log (N+1))^{2},\ \ \  (k\geq 1)
$$
where $\widetilde c$ depends only on $r,s$.\ \ \ \ \ $\Box$ 
\\
{\bf Proof of Corollary \ref{col1}.} According to Kolomogorov's theorem 
(see, e.g., [Gn], Ch. VI, Sec. 34) the condition
$$ 
\sum_{i=1}^{\infty}\frac{D(N_{i})}{i^{2}}<\infty
$$
guarantees fulfilment of the strong law of large numbers for the sequence
$\{N_{k}\}$.
This condition is a direct consequence of Theorem \ref{te2}. Furthermore,
Theorem \ref{te2} implies
$$
\frac{1}{n}\sum_{k=1}^{n}E(N_{k})\leq c\log M\log(N+1).
$$
Then the required statement follows from the strong law of large numbers
and the above inequality. \ \ \ \ \ $\Box$ \\
\begin{Ex}\label{ex1'}
{\rm Let $f:\Di_{1}\times B_{c}(0,r)\longrightarrow B_{c}(0,K)$ be a 
holomorphic mapping, where $B_{c}(0,r)$ and $B_{c}(0,K)$ are complex balls 
in $\Co^{N}$. Let $Df_{z}(v)$ denote the Jacobi 
matrix of $f(z,\cdot):B_{c}(0,r)\longrightarrow B_{c}(0,K)$, 
$z\in\Di_{1}$, at $v\in B_{c}(0,r)$. Assume that
$$
K_{1}:=\sup_{z\in\Di_{1}}\sup_{v\in B_{c}(0,r)}||Df_{z}(v)||_{2}<\infty \ .
$$
Here $||A||_{2}$ stands for the norm of a linear mapping $A: l_{2}(\Co^{N})
\longrightarrow l_{2}(\Co^{N})$. \\
Let $f:=(f_{1},...,f_{N})$.
Consider the system of ordinary differential equations
$$
\frac{dx_{i}}{dz}=f_{i}(z,x_{1},...,x_{N}),\ \ \ \ \ i=1,...,N
$$
with an initial condition 
$$
(x_{1}(0),...,x_{N}(0))=v\in B_{c}(0,3r/4)\ .
$$
In the special case that the system of equations above expresses the
property that the ${\rm N^{th}}$ 
derivative of a function is equal to zero, the solutions
are polynomials of degree $\leq N$. Therefore the inequalities below
apply to describe distribution of zeros of random polynomials.

For $|z|$ being sufficiently small one can solve the above system by 
the method of 
successive approximations. Namely, we choose $x_{0}(z,v)=v$ and then
$$
x_{n+1}(z,v):=v+\int_{0}^{z}f(w,x_{n}(w,v))dw\ .
$$
Here the integral is taken over the segment connecting 0 and $z$.
Further, if we take $R:=R(r,K,K_{1})<min\{r/4K,1/K_{1},1\}$ 
then $\{x_{n}\}_{n\geq 0}$ 
converges uniformly on $\Di_{R}\times B_{c}(0,3r/4)$ to a holomorphic  
mapping $x:\Di_{R}\times B_{c}(0,3r/4)\longrightarrow B_{c}(0,r)$ which
solves the above initial value problem (and this solution is unique).\\
Let $x_{1}$ be the 
first coordinate of $x$. Denote by
${\cal N}_{x_{1}}(v)$ the number of zeros of 
$x_{1v}:=x_{1}(\cdot,v)$, $v\in B_{c}(0,r/2)$, in a closed disk
$\overline{\Di}_{t}$, $t<R$. 

By definition we have
$$ 
M:=\sup_{v\in B_{c}(0,3r/4)}\sup_{z\in\Di_{R}}|x_{1v}(z)|\leq r\ \ \
{\rm and}\ \ \  M_{1}:=\sup_{v\in B_{c}(0,r/2)}|x_{1v}(0)|=r/2
$$
such that $M/M_{1}\leq 2$.
Therefore according to Theorem B,}

for every $T\geq 0$
$$
|\{v\in B_{c}(0,r/2)\ ;\ {\cal N}_{x_{1}}(v)\geq T\}|\leq 
Nc_{1}e^{-c_{2}T}|B_{c}(0,r/2)|,
$$
where $c_{1}, c_{2}$ depend only on $t/R,r$.\\
{\rm Consider now a particular case of the scheme described in 
Corollary \ref{col1}. Let $f_{v}^{(k)}(z)=x_{1v}(z),\ k=1,2,...$.
We define distributions of $N_{k}$ as follows:
for a nonnegative $l$ probability
$$
P(N_{k}=l):=\frac{1}{|B_{c}(0,r/2)|}|\{v\in B_{c}(0,r/2)\ ;\ 
{\cal N}_{x_{1}}(v)=l\}|\ .
$$ 
Then Corollary \ref{col1} implies that the inequality}
$$
\limsup_{n\to\infty}\frac{1}{n}\sum_{k=1}^{n}N_{k}\leq c\log(N+1)
$$
{\rm holds with probability one. Here $c$ depends on $t/R$ and $r$ only.\\
The constant on the right estimates the expected number of zeros in
$\overline{\Di}_{t}$ of a random function $x_{1v},\ v\in B_{c}(0,r/2)$.
In the special case that our system of ODE's depends linearly on the
parameter $v$, e.g. in the polynomial case mentioned above, one can obtain
sharper estimates by making use of explicit formulae (cf. [EK]).
}
\end{Ex}
{\bf Proof of Theorem \ref{central}.}
A simple calculation based on the inequality of Theorem B shows
that the sequence $\{E(|N_{k}-E(N_{k})|^{3})\}$ of absolute moments of third 
order is bounded.
We will show now that

{\em there is an} $\epsilon>0$ 
{\em such that} $D(N_{k})\geq\epsilon$ {\em for every} $k$.\\
This and the above boundedness of
third moments 
imply Lindeberg's condition  for the sequence $\{N_{k}-E(N_{k})\}$ 
(see, e.g., [Gn], Ch. VIII, Sec. 42). Therefore the central limit 
theorem will be valid in this case.

Assume, to the contrary, that there exists a subsequence 
$\{N_{k_{i}}\}_{i\geq 1}$ such that $\lim_{i\to\infty}D(N_{k_{i}})=0$.
Since ${\cal H}(M,r,s)$ is a compact in the topology of 
uniform convergence on compact subsets of $B_{c}(0,r)\times\Di_{1}$,
we can assume, without loss of generality, that $f^{(k_{i})}$ converges
in this topology to a function $f\in {\cal H}(M,r,s)$. 
Similarly to the definition of
${\cal N}_{f^{(k)}}$, introduce the function 
${\cal N}_{f}$ counting the number of 
zeros of $f_{v}$ $(v\in B_{c}(0,1))$ in $\overline{\Di}_{s}$.
We will show that under the above assumption the following result holds.
\begin{Lm}\label{N}
${\cal N}_{f}$ equals a constant almost everywhere on $B_{c}(0,1)$.
\end{Lm}

Based on this lemma let us, first, complete the proof of the theorem.
 
According to Proposition \ref{pr1} ${\cal N}_{f}$ is a
nonnegative upper semicontinuous on $B_{c}(0,r)$ function assigning only 
integer values. This and assumption $(a)$ of the theorem imply that the 
set $\{v\in \overline{B}_{c}(0,1)\ ; {\cal N}_{f}(v)=0\}$ is nonempty 
relatively open in $\overline{B}_{c}(0,1)$. Thus ${\cal N}_{f}=0$ almost 
everywhere on 
$B_{c}(0,1)$. Further, for some $s''$ satisfying $s'<s''<s$ consider 
the function $\widetilde {\cal N}_{f}$ which counts the number of zeros of 
$f_{v}$ in the open disk $\Di_{s''}$.  Then $\widetilde {\cal N}_{f}$ is 
bounded from above and lower semicontinuous on 
$\overline{B}_{c}(0,1)\setminus V_{f}$ (see Remarks 
\ref{fef} and \ref{open}); here $V_{f}:=\{v\in B_{c}(0,r)\ ; 
f_{v}=0\}$ is a proper analytic subset of $B_{c}(0,r)$.  Clearly, 
$0\leq \widetilde {\cal N}_{f}\leq {\cal N}_{f}$ which implies 
$\widetilde {\cal N}_{f}=0$ 
almost everywhere on $\overline{B}_{c}(0,1)$.  Observe now that  
condition $(b)$ of the theorem implies existence of $\widetilde 
v\in\overline{B}_{c}(0,1)$ such that the holomorphic function 
$f_{\widetilde v}\neq 0$ on $\Di_{s''}$ and has at 
least one zero there. In particular, $\widetilde v\not\in V_{f}$ and 
$\widetilde {\cal N}_{f}(\widetilde v)\geq 1$. The latter shows that 
$\widetilde {\cal N}_{f}$ attains its maximum (which is $>0$) on a nonempty 
relatively open subset of $\overline{B}_{c}(0,1)\setminus V_{f}$.  It 
contradicts to the assumption $\widetilde {\cal N}_{f}=0$ 
almost everywhere on $B_{c}(0,1)$. This proves that   
$\inf_{k}D(N_{k})=\epsilon>0$. 

The theorem is proved.\ \ \ \ \ $\Box$\\
{\bf Proof of Lemma \ref{N}.}
Let us consider the integral
\begin{equation}\label{integ}
\frac{1}{2\pi i}\int_{\Gamma}\frac{\frac{d}{dz}f_{v}(z)}{f_{v}(z)}dz;
\end{equation}
where $\Gamma:=\partial\overline{\Di}_{s}$.
Note that the expression above makes sense only for 
$v\in B_{c}(0,1)$ satisfying $\min_{z\in\Gamma}|f_{v}(z)|\neq 
0$. Denote the set of these points by $R$.
If $v\in R$ then (\ref{integ}) coincides with ${\cal N}_{f}(v)$. Consider
now the set $V\subset B_{c}(0,1)\times\Gamma$ defined by
$$
\{(v,z)\in B_{c}(0,1)\times\Gamma \ ; f_{v}(z)=0\}.
$$
Clearly, $V$ is a real analytic subset of $B_{c}(0,1)\times\Gamma$.
Moreover, real dimension of $V$ is at most $2N-1$ (recall that
$dim_{\Co}B_{c}(0,1)=N$). In fact, consider the natural projection
$\pi_{2}:V\longrightarrow\Gamma$. By the definition of $V$, for any
$\gamma\in\Gamma$ the set
$\pi_{2}^{-1}(\gamma)$ is a complex analytic subset of 
$B_{c}(0,1)$. From assumption $(a)$ of the theorem it follows that 
$\pi_{2}^{-1}(\gamma)$ does not coincide with $B_{c}(0,1)$. In 
particular, $dim_{\Re}\pi_{2}^{-1}(\gamma)\leq 2N-2$. Applying 
Sard's theorem to $\pi_{2}$ we then conclude that $dim_{\Re}V\leq
dim_{\Re}\Gamma + 2N-2=2N-1$. Further, the natural projection
$\pi_{1}:V\longrightarrow B_{c}(0,1)$ maps $V$ onto a subanalytic 
subset $\pi_{1}(V)$ of $B_{c}(0,1)$ and therefore $dim_{\Re}\pi_{1}(V)\leq
dim_{\Re}V\leq 2N-1$. In particular, we obtain that 
$R(=B_{c}(0,1)\setminus\pi_{1}(V))$ is everywhere dense in $B_{c}(0,1)$.
In addition,  the sequence
$\{f^{(k_{i})}_{v}\}_{i\geq 1}$ converges to $f_{v}$ uniformly on 
$\Gamma$ for every $v\in R$.  Hence, there is a number 
$i_{0}=i_{0}(v)$ such that for any
$i\geq i_{0}$ the function $f_{v}^{(k_{i})}$ has no zeros on $\Gamma$.  
From here and (\ref{integ}) it follows that ${\cal N}_{f^{(k_{i})}}$ converges to
${\cal N}_{f}$ almost everywhere on $B_{c}(0,1)$. 
Moreover, from Theorem \ref{te2} it follows that the sequence 
$\{E(N_{k})\}$ of nonnegative numbers is uniformly bounded from above.
Thus, without loss of generality, we can assume that 
$\lim_{i\to\infty}E(N_{k_{i}})=c$ for some nonnegative $c$.
This implies that the sequence of functions 
$\{({\cal N}_{f^{(k_{i})}}-E(N_{k_{i}}))^{2}\}_{i\geq 1}$ converges to
$({\cal N}_{f}-c)^{2}$ almost everywhere on $B_{c}(0,1)$. From here and
Fatou's lemma one obtains
$$
\int_{B_{c}(0,1)}({\cal N}_{f}-c)^{2}dx=
\int_{B_{c}(0,1)}\lim_{i\to\infty}
({\cal N}_{f^{(k_{i})}}-E(N_{k_{i}}))^{2}dx\leq |B_{c}(0,1)|\liminf_{i\to\infty}
D(N_{k_{i}})=0.
$$
This implies ${\cal N}_{f}=c$ almost everywhere on $B_{c}(0,1)$. \ \ \ \ \ 
$\Box$\\ 
{\bf Proof of Theorem \ref{polin}.} Consider, first, 
$\overline{\Di}_{1-\epsilon}=A_{\epsilon}\cap\Di_{1}$ and estimate 
expectation $E_{1,k}$ of the number of zeros of $P_{k,v}$ in 
$\overline{\Di}_{1-\epsilon}$. To this end we apply Bernstein's 
doubling inequality for polynomials along with the assumptions of the 
theorem to obtain
$$
\sup_{||v||\leq 2}|a_{ik}(v)|\leq 2^{d(k)}\sup_{||v||\leq 1}|a_{ik}(v)|
\leq 2^{d(k)},\ \ \ \ (1\leq 0\leq k).
$$
Then
$$
\sup_{z\in\Di_{1-\epsilon/2}}\sup_{||v||\leq 2}|P_{k,v}(z)|\leq
\max_{i}\sup_{||v||\leq 2}|a_{ik}(v)|\sum_{s=0}^{k}(1-\epsilon/2)^{s}<
\frac{2}{\epsilon}2^{d(k)}.
$$
Similarly, according to assumption $(c)$ of the theorem 
$$
\sup_{z\in\Di_{1-\epsilon}}\sup_{||v||\leq 1}|P_{k,v}(z)|\geq
\sup_{||v||\leq 1}|P_{k,v}(0)|=\sup_{||v||\leq 1}|a_{kk}(v)|=1.
$$
These inequalities show that, up to some dilation of the space of 
parameters and $\Co$ with coefficients depending only on $\epsilon$, 
the family $\{P_{k,v}\}_{v}$ belongs to ${\cal H}(M,r,s)$ with
$M=\frac{2^{d(k)+1}}{\epsilon}$ and with corresponding $r,s$ depending
on $\epsilon$. Applying now Theorem
\ref{te2} we estimate  
the expectation $E_{1,k}$ of the number of zeros of $P_{k,v}$ in 
$\overline{\Di}_{1-\epsilon}$ by $c(\epsilon)d(k)\log(n(k)+1)$.\\
To estimate expectation $E_{2,k}$
of the number of zeros of $P_{k,v}$ outside of
$\Di_{1+\epsilon}$ consider the family of polynomials
$P_{k,v}'(z):=z^{k}P_{k,v}(1/z)$. This, clearly, reduces the problem
to that of estimating the expectation of the number of zeros of 
$P_{k,v}'$ in $\overline{\Di}_{1-\epsilon}$. As above we get in this
case the inequality
$E_{2,k}\leq c(\epsilon)d(k)\log(n(k)+1)$.\\
Finally, the required expectation $E(\widetilde N_{k})= 
k-E_{1,k}-E_{2,k}\geq k(1-2c(\epsilon)\frac{d(k)\log(n(k)+1)}{k})$.
It remains to note that from here and assumption $(a)$ of the theorem it 
follows 
$E(\widetilde N_{k})=k(1-o_{\epsilon}(1))$ as $k\to\infty$.

The theorem is proved.\ \ \ \ \ $\Box$ 

I would like to thank Professors Yu.Brudnyi,
M.Goldstein, P.Milman and Y.Yomdin 
for useful discussions.
 
Deparment of Mathematics, Ben-Gurion University of the Negev,\\
Beer-Sheva 84105, Israel;\ \ \ \ \ \ \ \ email: brudnyi@cs.bgu.ac.il
\end{document}